\documentclass[11pt]{amsart}
\usepackage{amsmath, amssymb}

\newcommand{\RR}{\ensuremath{\mathbb{R}}}

\newcommand{\hgot}{\ensuremath{\mathfrak{h}}}

\newcommand{\ggot}{\ensuremath{\mathfrak{g}}}
\newcommand{\mgot}{\ensuremath{\mathfrak{m}}}
\newcommand{\xgot}{\ensuremath{\mathfrak{X}}}
\newcommand{\ogot}{\ensuremath{\mathfrak{o}}}
\newcommand{\sgot}{\ensuremath{\mathfrak{s}}}

\DeclareMathOperator{\Ker}{Ker} \DeclareMathOperator{\rank}{rank}
\DeclareMathOperator{\dm}{dim} \DeclareMathOperator{\dg}{diag}
\DeclareMathOperator{\dt}{det} \DeclareMathOperator{\ident}{id}
 \DeclareMathOperator{\Tr}{Tr}

\begin{document}

\author[V. V. Balashchenko, A. Sakovich]{Vitaly V. Balashchenko, Anna Sakovich}

\title[Invariant $f$-structures
on the flag manifolds]{Invariant $f$-structures on the flag
manifolds {\boldmath $SO(n)/SO(2) \times SO(n-3)$}}

\address{Vitaly V. Balashchenko\newline Faculty of Mathematics and Mechanics\newline
Belarusian State University\newline F.Scorina av.~4\newline Minsk
220050, BELARUS }

\email{balashchenko@bsu.by; vitbal@tut.by}

\address{Anna Sakovich\newline Faculty of Mathematics and Mechanics\newline
Belarusian State University\newline F.Scorina av.~4\newline Minsk
220050, BELARUS }

\email{anya\_sakovich@tut.by}



\begin{abstract}
We consider manifolds of oriented flags $SO(n)/SO(2) \times
SO(n-3)$ $(n \geq 4)$ as $4$- and $6$-symmetric spaces and
indicate characteristic conditions for invariant Riemannian
metrics under which the canonical $f$-structures on these
homogeneous $\Phi$-spaces belong to the classes ${\bf Kill\, f}$,
${\bf NKf}$, and ${\bf G_1 f}$  of generalized Hermitian geometry.
\end{abstract}

\maketitle

\section{Introduction}

An important place among homogeneous manifolds is occupied by
homogeneous $\Phi$-spaces \cite{F,BS2} of order $k$ (which are
also referred to as $k$-symmetric spaces \cite{Ko}), i.e. the
homogeneous spaces generated by Lie group automorphisms $\Phi$
such that $\Phi ^k=\ident$. Each $k$-symmetric space has an
associated object, the commutative algebra $\mathcal{A}(\theta)$
of canonical affinor structures \cite{BS1,BS2}. In its turn,
$\mathcal{A}(\theta)$ contains well-known classical structures, in
particular, $f$-structures in the sense of K.Yano \cite{Y}. It
should be mentioned that an $f$-structure compatible with a
(pseudo-)Riemannian metric is known to be one of the central
objects in the concept of generalized Hermitian geometry
\cite{Ki}.

From this point of view it is interesting to consider manifolds of
oriented flags of the form
\begin{equation}\label{space}
SO(n)/SO(2) \times SO(n-3) \;\;\; (n \geq 4)
\end{equation}
as they can be generated by automorphisms of any even finite order
$k \geq 4$. At the same time, it can be proved that an arbitrary
invariant Riemannian metric on these manifolds is (up to a
positive coefficient) completely determined by the pair of
positive numbers $(s,t)$. Therefore, it is natural to try to find
characteristic conditions imposed on $s$ and $t$ under which
canonical $f$-structures on homogeneous manifolds (\ref{space})
belong to the main classes of $f$-structures in the generalized
Hermitian geometry. This question is partly considered in the
paper.

The paper is organized as follows.

In Section 2,  basic notions and results related to homogeneous
regular $\Phi$-spaces and canonical affinor structures on them are
collected. In particular, this section includes a precise
description of all canonical $f$-structures on homogeneous
$k$-symmetric spaces.

In Section 3, we dwell on the main concepts of generalized
Hermitian geometry and consider the special classes of metric
$f$-structures such as ${\bf Kill \, f}$, ${\bf NKf}$, and ${\bf
G_1 f}$.

In Section 4, we describe manifolds of oriented flags of the form
$$SO(n)/\underbrace{SO(2)\times \dotsb \times SO(2)}_m \times
SO(n-2m-1)$$ and construct inner automorphisms by which they can
be generated.

In Section 5, we describe the action of the canonical
$f$-structures on the flag manifolds of the form (\ref{space})
considered as homogeneous $\Phi$-spaces of orders 4 and 6.

Finally, in Section 6, we indicate characteristic conditions for
invariant Riemannian metrics on the flag manifolds (\ref{space})
under which the canonical $f$-structures on these homogeneous
$\Phi$-spaces belong to the classes ${\bf Kill \, f}$, ${\bf
NKf}$, and ${\bf G_1 f}$.

\section{ Canonical structures on regular
$\Phi$-spaces}

We start with some basic definitions and results related to
homogeneous regular $\Phi$-spaces and canonical affinor
structures. More detailed information can be found in \cite{S1},
\cite{F}, \cite{Ko}, \cite{BS2}, \cite{B2} and some others.

Let $G$ be a connected Lie group, $\Phi$ its automorphism.  Denote
by $G^{\Phi}$ the subgroup consisting of all fixed points of
$\Phi$ and by $G^{\Phi}_0$ the identity component of $G^{\Phi}$.
Suppose a closed subgroup $H$ of $G$ satisfies the condition
$$G^{\Phi}_0 \subset H \subset G^{\Phi}.$$ Then $G/H$ is called a
{\it homogeneous $\Phi$-space} \cite{F,BS2}.

Among homogeneous $\Phi$-spaces a fundamental role is played by
{\it homogeneous $\Phi$-spaces of order $k$} ($\Phi ^k=\ident$)
or, in the other terminology, {\it homogeneous $k$-symmetric
spaces} (see \cite{Ko}).

Note that there exist homogeneous $\Phi$-spaces that are not
reductive. That is why so-called regular $\Phi$-spaces first
introduced by N.A.Stepanov \cite{S1} are of fundamental
importance.

Let $G/H$ be a homogeneous $\Phi$-space, $\mathfrak{g}$ and
$\mathfrak{h}$ the corresponding Lie algebras for $G$ and $H$,
$\varphi=d{\Phi}_e$ the automorphism of $\mathfrak{g}$. Consider
the linear operator $A=\varphi-\ident$ and the Fitting
decomposition $\mathfrak{g}=\mathfrak{g}_0\oplus\mathfrak{g}_1$
with respect to $A$, where   $\mathfrak{g}_0$ and $\mathfrak{g}_1$
denote $0$- and $1$-component of the decomposition respectively.
Further, let $\varphi=\varphi_{s}\:\varphi_{u}$ be the Jordan
decomposition, where $\varphi_{s}$ and $\varphi_{u}$ is a
semisimple and unipotent component of $\varphi$ respectively,
$\varphi_{s}\:\varphi_{u}=\varphi_{u}\:\varphi_{s}$. Denote by
$\mathfrak{g}^{\gamma}$ a subspace of all fixed points for a
linear endomorphism $\gamma$ in $\mathfrak{g}$. It is clear that
$\mathfrak{h}=\mathfrak{g}^{\varphi}=Ker\,A$,
$\mathfrak{h}\subset\mathfrak{g}_0$,
$\mathfrak{h}\subset\mathfrak{g}^{\varphi_s}$.

{\bf Definition\;1}\;\cite{F,S1,BS2,B2}.\; A homogeneous
$\Phi$-space $G/H$ is called a {\it regular $\Phi$-space} if one
of the following equivalent conditions is satisfied:
\begin{enumerate}
\item $\mathfrak{h}=\mathfrak{g}_0$. \item
$\mathfrak{g}=\mathfrak{h}\oplus{A}\mathfrak{g}$. \item The
restriction of the operator $A$ to ${A}\mathfrak{g}$ is
non-singular. \item $A^2X=0\Longrightarrow{A}X=0$ for all
$X\in\mathfrak{g}$. \item The matrix of the automorphism $\varphi$
can be represented in the form $\left(\begin{array}{cc}
  E & 0 \\
  0 & B
\end{array}\right),$ where the matrix $B$ does not admit the eigenvalue
$1$. \item $\mathfrak{h}=\mathfrak{g}^{\varphi_s}$.
\end{enumerate}

A distinguishing feature of a regular $\Phi$-space $G/H$ is that
each such space is reductive, its reductive decomposition being
$\ggot=\hgot \oplus A \ggot$ (see \cite{S1}). $\ggot=\hgot \oplus
A \ggot$ is commonly referred to as the {\it canonical reductive
decomposition} corresponding to a regular $\Phi$-space $G/H$ and
$\mgot=A \ggot$ is the {\it canonical reductive complement}.

It should be mentioned that any homogeneous $\Phi$-space $G/H$ of
order $k$ is regular (see \cite{S1}), and, in particular, any
$k$-symmetric space is reductive.

Let us now turn to canonical $f$-structures on regular
$\Phi$-spaces.

An {\it affinor structure} on a smooth manifold is a tensor field
of type $(1,1)$ realized as a field of endomorphisms acting on its
tangent bundle. It is known that any invariant affinor structure
$F$ on a homogeneous manifold $G/H$ is completely determined by
its value $F_o$ at the point $o=H$, where $F_o$ is invariant with
respect to $Ad(H)$. For simplicity, further we will not
distinguish an invariant structure on $G/H$ and its value at $o=H$
throughout the rest of the paper.

Let us denote by $\theta$ the restriction of $\varphi$ to
$\mgot$.\vspace{6pt}

{\bf Definition\;2}\;\cite{BS1,BS2}.\;An invariant affinor
structure $F$ on a regular $\Phi$-space $G/H$ is called {\it
canonical} if its value at the point $o=H$ is a polynomial in
$\theta$.\vspace{6pt}

Remark that the set $\mathcal{A}(\theta)$ of all canonical
structures on a regular $\Phi$-space $G/H$ is a commutative
subalgebra of the algebra $\mathcal{A}$ of all invariant affinor
structures on $G/H$. This subalgebra contains well-known classical
structures such as {\it almost product structures} ($P^2=\ident$),
{\it almost complex structures} ($J^2=-\ident$), {\it
$f$-structures} ($f^3+f=0$).

The sets of all canonical structures of the above types were
completely described in \cite{BS1} and \cite{BS2}. In particular,
for homogeneous $k$-symmetric spaces the precise computational
formulae were indicated. For future reference we cite here the
result pertinent to $f$-structures and almost product structures
only. Put
$$u=\left\{\begin{array}{ll}
 n \; &\text{if } \; k=2n+1, \\
 n-1\; & \text{if } \; k=2n.
 \end{array}\right.$$

{\bf Theorem 1}\;\cite{BS1,BS2}.\;{\it Let $G/H$ be a
homoge\-neous $\Phi$-space of order $k$ $(k \geq 3)$.}
\begin{itemize}
\item[1)] {\it All non-trivial canonical $f$-structures on $G/H$ can be given
by the operators
\begin{equation*}
f(\theta)=\frac{2}{k}\sum _{m=1} ^{u} \left(\sum _{j=1} ^{u} \zeta
_j \sin \frac{2 \pi m j}{k}\right)(\theta ^m-\theta^{k-m}),
\end{equation*}
where $\zeta_j \in \{1,0,-1\},\; j=1,2, \dotsc, u$, and not all
$\zeta_j$ are equal to zero. }
\item[2)]{\it All
canonical almost product structures $P$ on $G/H$ can be given by
polynomials $P(\theta)=\sum_{m=0} ^{k-1} a_m \theta ^m$, where:}
\begin{itemize}
\item [a)] {\it if $k=2n+1$, then $$a_m=a_{k-m}=\frac{2}{k}\sum _{j=1}^u
\xi _j \cos \frac{2\pi m j}{k};$$}

\item [b)] {\it if $k=2n$, then
$$a_m=a_{k-m}=\frac{1}{k}\left(2\sum _{j=1}^u \xi _j \cos
\frac{2\pi m j}{k}+(-1)^{m}\xi _n \right).$$}
\end{itemize}
{\it Here the numbers $\xi _j$, $j=1,2,\dotsc ,u$, take their
values from the set $\{-1,1\}$.}
\end{itemize}

The results mentioned above were particularized for homogeneous
$\Phi$-spaces of smaller orders $3$, $4$, and $5$ (see
\cite{BS1,BS2}). Note that there are no fundamental obstructions
to considering of higher orders $k$. Specifically, for future
consideration we need the description of canonical $f$-structures
and almost product structures on homogeneous $\Phi$-spaces of
orders $4$ and $6$ only.

{\bf Corollary 1}\;\cite{BS1,BS2}.\;{\it Any homogeneous
$\Phi$-space of order $4$ admits (up to sign) the only canonical
$f$-structure $$f_0(\theta)=\frac{1}{2}(\theta-\theta^3)$$ and the
only almost product structure $$P_0(\theta)=\theta ^2.$$}

{\bf Corollary 2}.\;{\it On any homogeneous $\Phi$-space of order
$6$ there exist (up to sign) only the following canonical
$f$-structures:

\begin{align*} f_1(\theta)&
=\frac{1}{\sqrt{3}}(\theta-\theta^5), & f_2(\theta)&
=\frac{1}{2\sqrt{3}}(\theta-\theta^2+\theta^4-\theta^5),\\
f_3(\theta)&
=\frac{1}{2\sqrt{3}}(\theta+\theta^2-\theta^4-\theta^5),&
f_4(\theta)& =\frac{1}{\sqrt{3}}(\theta^2-\theta^4)
\end{align*}
and only the following almost product structures:
\begin{align*}
P_1(\theta)& =-\ident,& P_2(\theta)&
=\frac{\theta}{3}+\theta^2+\frac{\theta^3}{3}+\theta
^4+\frac{\theta^5}{3},\\ P_3(\theta)& =\theta ^3, & P_4(\theta)&
=-\frac{2 \theta ^2}{3}+\frac{\theta ^3}{3}-\frac{2 \theta ^5}{3}.
\end{align*}}

\section{ Some important classes in generalized Hermitian
geometry}

The concept of generalized Hermitian geometry created in the 1980s
(see \cite{Ki}) is a natural consequence of the development of
Hermitian geometry. One of its central objects is a {\it metric
$f$-structure}, i.e. an $f$-structure compatible with a
(pseudo-)Riemannian metric $g=\langle \cdot, \cdot \rangle$ in the
following sense: $$\langle fX, Y \rangle + \langle X, fY \rangle
=0 \text{ for any } X,\; Y \in \xgot (M).$$ Evidently, this
concept is a generalization of one of the fundamental notions in
Hermitian geometry, namely, almost Hermitian structure $J$. It is
also worth noticing that the main classes of generalized Hermitian
geometry (see \cite{Ki,Gr1, Gr2, B2, B6}) in the special case
$f=J$ coincide with those of Hermitian geometry (see \cite{GH}).

In what follows, we will mainly concentrate on the classes ${\bf
Kill \, f}$, ${\bf NKf}$, and ${\bf G_1 f}$ of metric
$f$-structures defined below.

A fundamental role in generalized Hermitian geometry is played by
a tensor $T$ of type (2,1) which is called a {\it composition
tensor} \cite{Ki}. In \cite{Ki} it was also shown that such a
tensor exists on any metric $f$-manifold and it is possible to
evaluate it explicitly:
$$T(X,Y)=\frac{1}{4}f(\nabla _{fX} (f)
fY-\nabla _{f^2 X} (f) f^2 Y), $$ where $\nabla$ is the
Levi-Civita connection of a (pseudo-)Riemannian manifold $(M,g)$,
$X,\,Y \in \xgot (M)$.

The structure of a so-called {\it adjoint $Q$-algebra} (see
\cite{Ki}) on $\xgot (M)$ can be defined by the formula
$X*Y=T(X,Y)$. It gives the opportunity to introduce some classes
of metric $f$-structures in terms of natural properties of the
adjoint $Q$-algebra. For example, if $T(X,X)=0$ (i.e. $\xgot (M)$
is an anticommutative $Q$-algebra) then $f$ is referred to as a
{\it $G_1 f$-structure}. ${\bf G_1 f}$ stands for the class of
$G_1 f$-structures.

A metric $f$-structure on $(M,g)$ is said to be a {\it Killing
$f$-structure} if $$\nabla _X (f) X=0 \text{ for any } X \in \xgot
(M)$$ (i.e. $f$ is a Killing tensor) (see \cite{Gr1, Gr2}). The
class of Killing $f$-structures is denoted by ${\bf Kill\, f}$.
The defining property of {\it nearly K\"ahler $f$-structures} (or
{\it $NKf$-structures}) is $$\nabla _{fX} (f) {fX}=0.$$ This class
of metric $f$-structures, which is denoted by ${\bf NKf}$, was
determined in \cite{B6} (see also \cite{B4,B3}). It is easy to see
that for $f=J$ the classes ${\bf Kill\, f}$ and ${\bf NKf}$
coincide with the well-known class ${\bf NK}$ of {\it nearly
K\"ahler structures} \cite{G1}.

The following relations between the classes mentioned are evident:
$${\bf Kill \, f} \subset  {\bf NKf} \subset {\bf G_1 f}.$$

A special attention should be paid to the particular case of
naturally reductive spaces. Recall that a homogeneous Riemannian
manifold $(G/H,g)$ is known to be a {\it naturally reductive
space} \cite{KN} with respect to the reductive decomposition
$\ggot=\hgot \oplus \mgot$ if
$$g([X,Y]_{\mgot},Z)=g(X,[Y,Z]_{\mgot}) \text{ for any } X,Y,Z \in
\mgot. $$ It should be mentioned that if $G/H$ is  a regular
$\Phi$-space, $G$ a semisimple Lie group then $G/H$ is a naturally
reductive space with respect to the (pseudo-)Riemannian metric $g$
induced by the Killing form of the Lie algebra $\ggot$ (see
\cite{S1}). In \cite{B4}, \cite{B3}, \cite{B5} and \cite{B6} a
number of results helpful in checking whether the particular
$f$-structure on a naturally reductive space belongs to the main
classes of generalized Hermitian geometry was obtained.

\section{ Manifolds of oriented flags}

In linear algebra  a {\it flag} is defined as a finite sequence
$L_0,\dotsc,L_n$ of subspaces of a vector space $L$ such that
\begin{equation}\label{flag}
L_0 \subset L_1 \subset \dotsb \subset L_n,
\end{equation}
$L_i \neq L_{i+1}, \; i=0,\dotsc,n-1$ (see \cite{KM}).

A flag (\ref{flag}) is known to be {\it full}\; if for any
$i=0,\dotsc ,n-1$ $\dm L_{i+1}=\dm L_i +1$. It is readily seen
that having fixed any basis $\{ e_1\dotsc ,e_n \}$ of $L$ we can
construct a full flag by setting $L_0=\{0\}, \;
L_i=\mathcal{L}(e_1,\dotsc ,e_i), \; i=1,\dotsc,n.$

We call a flag $L_{i_1} \subset L_{i_2} \subset \dotsb \subset
L_{i_n}$ (here and below the subscript denotes the dimension of
the subspace) {\it oriented} if for any $L_{i_j}$ and its two
basises $\{e_1,\dotsc,e_{i_j}\}$ and $\{e'_1,\dotsc,e'_{i_j}\}$
$\dt A
> 0$, where $e'_t=A e_t$ for any $t=1,\dotsc,i_j$. Moreover,
for any two subspaces $L_{i_k} \subset L_{i_j}$ their orientations
should be set in accordance.

{\bf Proposition 1.} {\it The set of all oriented flags $$L_1
\subset L_3 \subset \dotsb \subset L_{2m+1} \subset L_n=L$$ of a
vector space $L$ with respect to the action of $SO(n)$ is
isomorphic to $$SO(n)/\underbrace{SO(2)\times \dotsb \times
SO(2)}_m \times SO(n-2m-1).$$}

{\bf Proof.} Fix some basis $\{e_1,\dotsc,e_n \}$ in $L_n$.
Consider the isotropy subgroup $I_o$ at the point $$
o=(\mathcal{L}(e_1)\subset \mathcal{L}(e_1,e_2,e_3)\subset \dotsb
\subset \mathcal{L}(e_1,\dotsc,e_{2m+1}) \subset \mathcal{L}
(e_1,\dotsc ,e_n)). $$

By the definition for any $A \in I_o$ $$A:\mathcal{L}(e_1)
\rightarrow \mathcal{L}(e_1),$$ $$A:\mathcal{L}(e_1,e_2,e_3)
\rightarrow \mathcal{L}(e_1,e_2,e_3),\dotsc,$$
$$A:\mathcal{L}(e_1,\dotsc,e_{2m+1}) \rightarrow
\mathcal{L}(e_1,\dotsc,e_{2m+1}),$$ $$A:\mathcal{L}
(e_1,\dotsc,e_n) \rightarrow \mathcal{L} (e_1,\dotsc,e_n).$$

As $\{e_1,\dotsc,e_n \}$ is a basis, it immediately follows that
$$A:\mathcal{L}(e_1) \rightarrow \mathcal{L}(e_1),$$
$$A:\mathcal{L}(e_2,e_3) \rightarrow
\mathcal{L}(e_2,e_3),\dotsc,$$ $$A:\mathcal{L}(e_{2m},e_{2m+1})
\rightarrow \mathcal{L}(e_{2m},e_{2m+1}),$$ $$A:\mathcal{L}
(e_{2m+2},\dotsc,e_n) \rightarrow \mathcal{L}
(e_{2m+2},\dotsc,e_n).$$

Thus $L=L_n$ can be decomposed into the sum of $A$-invariant
subspaces $$L=\mathcal{L}(e_1) \oplus \mathcal{L}(e_2,e_3) \oplus
\dots \oplus \mathcal{L}(e_{2m},e_{2m+1})\oplus \mathcal{L}
(e_{2m+2},\dotsc,e_n).$$ The matrix of the operator $A$ in the
basis $\{e_1,\dotsc,e_n \}$ is cellwise-diagonal:$$A=\dg \{A^1
_{1\times 1} , A^3 _{2 \times 2},\dotsc,A^{2m+1} _{2 \times 2},
A^n _{(n-2m-1)\times (n-2m-1)} \}.$$ Since $A \in SO(n)$, its
cells $A^1 , A^3, \dotsc ,A^{2m+1} , A^n$ are orthogonal matrices.
All the flags we consider are oriented, thus for any $i \in
\{1,3,\dotsc,2m+1,n \}$  $\dt A^i >0$. This proves that $A^1=(1)$,
$A^3 \in SO(2),\dotsc,A^{2m+1} \in SO(2)$, $A^n \in SO(n-2m-1)$.

Therefore $I_o=\underbrace{SO(2) \times \dotsb \times SO(2)}_m
\times SO(n-2m-1).$ This completes the proof. \hfill $\Box$

{\bf Proposition 2}. {\it The manifold of oriented flags
$$SO(n)/\underbrace{SO(2)\times \dotsb \times SO(2)}_m \times
SO(n-2m-1)$$ is a homogeneous $\Phi$-space. It can be generated by
inner automorphisms $\Phi$ of any finite order $k$, where $k$ is
even, $k>2$ and $k\geq 2m-2$: $$\Phi:SO(n)\rightarrow SO(n), \; A
\rightarrow BAB^{-1}, \text{ where}$$ $$B=\dg \{1, \varepsilon
_1,\dotsc, \varepsilon_m, -1,\dotsc,-1 \}, $$ $$\varepsilon
_t=\begin{pmatrix}\cos \frac{2 \pi t}{k} & \sin \frac{2 \pi t}{k}
\\[4pt] -\sin \frac{2 \pi t}{k} & \cos \frac{2 \pi
t}{k}\end{pmatrix}.$$}

{\bf Proof.} Here $G=SO(n)$, $H=\underbrace{SO(2)\times \dotsb
\times SO(2)}_m \times SO(n-2m-1)$. We need to prove that the
group of all fixed points $G^{\Phi}$ satisfies the condition
$$G^{\Phi} _0 \subset H \subset G^{\Phi}.$$ By definition
$G^{\Phi}=\{A|B A B^{-1}=A\}=\{A|BA=AB\}$. Equating the
correspondent elements of $AB$ and $BA$ and solving systems of
linear equations it is possible to calculate that $$G^{\Phi}=\{\pm
1\} \times \underbrace{SO(2)\times \dots\times SO(2)}_m \times
SO(n-2m-1). $$ \hfill $\Box$

\section{ Canonical $f$-structures  on
$4$- and $6$-symmetric space $SO(n)/SO(2)\times SO(n-3)$}

Let us consider $SO(n)/SO(2) \times SO(n-3)$ $(n \geq 4)$ as a
homogeneous $\Phi$-space of order 4. According to Proposition 2 it
can be generated by the inner automorphism $\Phi: A \rightarrow B
A B^{-1},$ where $$B=\dg \left \{1,\left(\begin{array}{cc} 0 & 1
\\ - 1 & 0\end{array}\right),\underbrace{-1,\dotsc,
-1}_{n-3}\right \}.$$ Therefore (\ref{space}) is a reductive
space. It is not difficult to check that the canonical reductive
complement $\mgot$ consists of matrices of the form

$$S=\begin{pmatrix} 0 & s_{12} & s_{13} & s_{14} & \dots &
s_{1n}\\ -s_{12} & 0 & 0 & s_{24} & \dots & s_{2n}\\ -s_{13} & 0 &
0 & s_{34} & \dots & s_{3n}\\ -s_{14} & -s_{24} & -s_{34} & 0 &
\dots & 0\\  \hdotsfor[2]{6} \\ -s_{1n} & -s_{2n} & -s_{3n} & 0 &
\dots & 0\\
\end{pmatrix} \in \mgot.$$ According to Corollary 1 the only canonical
$f$-structure on this homogeneous $\Phi$-space
is determined by the formula
$$f_0(\theta)=\frac{1}{2}(\theta-\theta ^3).$$ Its action can be
written in the form: $$f_0:S \longrightarrow
\begin{pmatrix}0 & s_{13} & -s_{12} & 0 & \dots & 0\\
-s_{13} & 0 & 0 & -s_{34} & \dots & -s_{3n}\\ s_{12} & 0 & 0 &
s_{24} & \dots & s_{2n}\\ 0 & s_{34} & -s_{24} & 0 & \dots & 0\\
\hdotsfor[2]{6}\\ 0 & s_{3n} & -s_{2n} & 0 & \dots & 0\\
\end{pmatrix}.$$

Now let us consider (\ref{space}) as a 6-symmetric space generated
by the inner automorphism $\Phi: A \rightarrow B A B^{-1},$ where
$$B=\dg \left \{1,\left(\begin{array}{cc} \frac{1}{2} &
\frac{\sqrt{3}}{2}
\\[2pt] - \frac{\sqrt{3}}{2} & \frac{1}{2} \end{array}\right),\underbrace{-1,\dotsc
, -1}_{n-3}\right \}.$$ Taking Corollary 2 into account we can
represent the action of the canonical $f$-structures on this
homogeneous $\Phi$-space as follows:
$$f_1(\theta)=\frac{1}{\sqrt{3}}(\theta-\theta^5):S
\longrightarrow
\begin{pmatrix} 0 & s_{13} & -s_{12} & 0 & \dots & 0\\
-s_{13} & 0 & 0 & -s_{34} & \dots & -s_{3n}\\ s_{12} & 0 & 0 &
s_{24} & \dots & s_{2n}\\ 0 & s_{34} & -s_{24} & 0 & \dots &
0\\\hdotsfor[2]{6}\\ 0 & s_{3n} & -s_{2n} & 0 & \dots & 0\\
\end{pmatrix},$$

$$f_2(\theta)=\frac{1}{2\sqrt{3}}(\theta-\theta^2+\theta^4-\theta^5):S
\longrightarrow \begin{pmatrix} 0 & 0 & 0 & 0 & \dots & 0\\ 0 & 0
& 0 & -s_{34} & \dots & -s_{3n}\\ 0 & 0 & 0 & s_{24} & \dots &
s_{2n}\\ 0 & s_{34} & -s_{24} & 0 & \dots & 0\\ \hdotsfor[2]{6}\\
0 & s_{3n} & -s_{2n} & 0 & \dots & 0\\
\end{pmatrix},$$

$$f_3(\theta)=\frac{1}{2\sqrt{3}}(\theta+\theta^2-\theta^4-\theta^5):S
\longrightarrow \begin{pmatrix} 0 & s_{13} & -s_{12} & 0 & \dots &
0\\ -s_{13} & 0 & 0 & 0 & \dots & 0\\ s_{12} & 0 & 0 & 0 & \dots &
0\\ 0 & 0 & 0 & 0 & \dots & 0\\\hdotsfor[2]{6}\\ 0 & 0 & 0 & 0 &
\dots & 0\\
\end{pmatrix},$$

$$f_4(\theta)=\frac{1}{\sqrt{3}}(\theta^2-\theta^4):S
\longrightarrow
\begin{pmatrix} 0 & s_{13} & -s_{12} & 0 & \dots & 0\\
-s_{13} & 0 & 0 & s_{34} & \dots & s_{3n}\\ s_{12} & 0 & 0 &
-s_{24} & \dots & -s_{2n}\\ 0 & -s_{34} & s_{24} & 0 & \dots & 0\\
\hdotsfor[2]{6}\\ 0 & -s_{3n} & s_{2n} & 0 & \dots & 0\\
\end{pmatrix}.$$

\section{ Canonical $f$-structures and invariant Riemannian
metrics on $SO(n)/SO(2)\times SO(n-3)$}

Let us consider manifolds of oriented flags of the form
(\ref{space})  as $4$- and $6$-symmetric spaces. Our task is to
indicate characteristic conditions for invariant Riemannian
metrics under which the canonical $f$-structu\-res on these
homogeneous $\Phi$-spaces belong to the classes ${\bf Kill \, f}$,
${\bf NKf}$, and ${\bf G_1 f}$.

We begin with some preliminary considerations.

{\bf Proposition 3.} {\it The reductive complement $\mgot$ of the
homogeneous space $SO(n)/SO(2)\times SO(n-3)$ admits the
decomposition into the direct sum of $Ad(H)$-invariant irreducible
subspaces $\mgot=\mgot _1 \oplus \mgot _2 \oplus \mgot _3 $.}

{\bf Proof.} The explicit form of the reductive complement of
(\ref{space}) was indicated in Section 5. Put $$\mgot _1=\left
\{\left. \begin{pmatrix} 0 & a_{1} & a_{2} & 0 & \dots & 0\\
-a_{1} & 0 & 0 & 0 & \dots & 0\\ -a_2 & 0 & 0 & 0 & \dots & 0\\ 0
& 0 & 0 & 0 & \dots & 0\\ \hdotsfor[2]{6}\\ 0 & 0 & 0 & 0 & \dots
& 0\\
\end{pmatrix}\right| a_1,\; a_2 \in \RR \right \},$$

$$\mgot_2=\left \{\left. \begin{pmatrix} 0 & 0 & 0 & 0 & \dots &
0\\0 & 0 & 0 & c_1 & \dots & c_{n-3}\\ 0 & 0 & 0 & d_1 & \dots &
d_{n-3}\\ 0 & -c_1 & -d_1 & 0 & \dots & 0\\  \hdotsfor[2]{6}\\ 0 &
-c_{n-3} & -d_{n-3} & 0 & \dots & 0\\
\end{pmatrix}\right|\begin{array}{cc} c_1,\dotsc,c_{n-3} \in
\RR \\ d_1,\dotsc,d_{n-3} \in \RR \end{array} \right \},$$

$$\mgot_3=\left \{\left. \begin{pmatrix} 0 & 0 & 0 & b_1 & \dots &
b_{n-3}\\0 & 0 & 0 & 0 & \dots & 0\\ 0 & 0 & 0 & 0 & \dots & 0\\
-b_1 & 0 & 0 & 0 & \dots & 0\\  \hdotsfor[2]{6}\\ -b_{n-3} & 0 & 0
& 0 & \dots & 0\\
\end{pmatrix}\right| b_1,\dotsc,b_{n-3} \in \RR
 \right \}.$$

Since $SO(2) \times SO(n-3)$ is a connected Lie group, $\mgot _i$
($i=1,\,2,\,3$) is $Ad(H)$-invariant iff $[\hgot,\mgot_i ] \subset
\mgot_i$. It can easily be shown that this condition holds.

We claim that for any $i \in \{1,\,2,\,3\}$ there exist no such
non-trivial subspaces $\overline{\mgot}_i$ and $\widehat{\mgot}_i$
that $\mgot_i=\overline{\mgot}_i \oplus \widehat{\mgot}_i$ and
$[\hgot,\overline{\mgot}_i ] \subset \overline{\mgot}_i$,
$[\hgot,\widehat{\mgot}_i ] \subset \widehat{\mgot}_i$.

To prove this we identify $\mgot$ and $$
\{(a_{1},a_{2},b_{1},\dotsc,\,b_{n-3},c_{1},\dotsc,c_{n-3},d_{1},\dotsc,d_{n-3})\}.$$
In what follows we are going to represent any $H \in \hgot$ in the
form $$H=\dg \{0,\,H_1,\,H_2\}$$ where $$H_1=\begin{pmatrix} 0 & h
\\ -h & 0
\end{pmatrix},$$

\begin{equation}\label{H2}
H_2=\begin{pmatrix} 0 & h_{1 \, 2} & \dots & h_{1 \, n-3}
\\-h_{1 \, 2} & 0 & \dots & h_{2 \, n-3} \\ \hdotsfor[2]{4}
\\  -h_{1 \,n-3}& -h_{2 \, n-3} & \dots  & 0 \end{pmatrix}.
\end{equation}

Put $F(H)(M)=[H,M]$ for any $H \in \hgot$, $M \in \mgot$. In the
above notations we have

$$F(H)|_{\mgot _1}: (a_1 \;a_2)^{T} \rightarrow H_1 (a_1\;
a_2)^{T},$$

\begin{multline*}
F(H)|_{\mgot_2}: (c_1 \dotso c_{n-3} \; d_1 \dotso d_{n-3})^T
\rightarrow \\ \rightarrow \left( \begin{array}{cc} H_2 & h E
\\ -h E & H_2 \end{array} \right ) (c_1 \dotso c_{n-3} \; d_1 \dotso d_{n-3})^T,
\end{multline*}

$$F(H)|_{\mgot_3}:(b_1 \; \dots\; b_{n-3} )^T \rightarrow H_2 (b_1
\dotso b_{n-3} )^T.$$

First, let us prove that $\mgot_3$ cannot be decomposed into the
direct sum of $Ad(H)$-invariant subspaces.

The proof is by reductio ad absurdum. Suppose there exists an
$Ad(H)$-invariant subspace $W \subset \mgot _3$. This implies that
for any $H_2$ of the form (\ref{H2}) and $x=(x_1 \dotso x_{n-3})^T
\in W$ $H_2 x$ belongs to $W$.

It is possible to choose a vector $v_1=(\alpha _1 \dotso
\alpha_{n-3})^T \in W$ such that $\alpha _1 \neq 0$. Indeed, the
nonexistence of such a vector yields that for any $w=(w_1 \dotso
w_{n-3})^T \in W$ $w_1=0$. Take such $w \in W$ that for some $1< i
\leq n-3$ $w_i \neq 0$ and the skew-symmetric matrix $K=\{k_{i \,
j}\}$ with all elements except $k_{1\, i}=-k_{i \, 1}=1$ equal to
zero. Then $Kw=(w_i \; \ast \dotso \ast) \notin W.$

Consider the following system of vectors $\{ v_1, \dotsc ,v_{n-3}
\}$, where $$v_2=\begin{pmatrix}0 & 1 & 0 & \dots & 0\\ -1 & 0 & 0
& \dots & 0\\   0 & 0 & 0 & \dots & 0 \\ \hdotsfor[2]{5}
\\ 0 & 0 & 0 & \dots & 0
\end{pmatrix} v_1=(\alpha_2 \;- \alpha_1 \; 0 \dotso 0)^T,$$

$$v_3=\begin{pmatrix} 0 & 0 & 1 & \dots & 0\\ 0 &  0 & 0 & \dots &
0\\   -1 & 0 & 0 & \dots & 0 \\ \hdotsfor[2]{5}
\\ 0 & 0 & 0 & \dots & 0
\end{pmatrix} v_1=(\alpha_3 \; 0  \; -\alpha_1 \dotso 0)^T,\;\dots\;,$$

$$v_{n-3}=\begin{pmatrix} 0 & 0 & \dots & 0 & 1\\ 0 & 0 & \dots &
0 & 0\\  \hdotsfor[2]{5} \\ 0 & 0 & \dots & 0 & 0
\\ -1 & 0 & \dots & 0 & 0
\end{pmatrix} v_1=(\alpha_n \; 0 \dotso 0 \; -\alpha _1)^T.$$

Obviously, $\dm \mathcal{L} (v_1, \dots, v_{n-3})=$ $$=\rank
\begin{pmatrix} \alpha _1 & \alpha _2  & \alpha _3 &
\dots & \alpha _{n-3}\\ \alpha _2 & -\alpha _1 & 0 & \dots & 0\\
\alpha _3 & 0 & -\alpha_1 & \dots & 0
\\ \hdotsfor[2]{5}
\\ \alpha _{n-3} & 0 & 0 & \dots & -\alpha _1
\end{pmatrix}=n-3.$$ This contradicts our assumption.

Continuing the same line of reasoning, we see that neither $\mgot
_1$ nor $\mgot _2$ can be decomposed into the sum of
$Ad(H)$-invariant summands. \hfill $\Box$

It is not difficult to check that the space in question possesses
the following property.

{\bf Proposition\;4.}
\begin{equation}\label{locrel} [\mgot _i
,\mgot _{i+1}] \subset \mgot _{i+2}\; ( \text{{\it modulo 3}}).
\end{equation}

Denote by $g_0$ the naturally reductive metric generated by the
Killing form $B$: $g_0=-B|_{\mgot \times \mgot}$. In our case
$B=-(n-1)\Tr X^T Y$, $X,Y \in \sgot \ogot (n)$.

{\bf Proposition\;5.}\;{\it The decomposition $\hgot \oplus \mgot
_1 \oplus \mgot _2 \oplus \mgot _3$ is $B$-orthogonal.}

{\bf Proof.} For the explicit form of $\mgot$ and $\hgot$ see
Section 5 and Section 6. It can easily be seen that for any  $X
\in \mgot$, $Y \in \hgot$ $\Tr X^T Y=0$. It should also be noted
that it was proved in \cite{S1} that $\hgot$ is orthogonal to
$\mgot$ with respect to $B$.

For any almost product structure $P$ put $$\mgot^-=\{X \in
\mgot\;| P(X)=-X\}, \; \mgot ^+=\{X \in \mgot\;|P(X)=X\}.$$
Suppose that $P$ is compatible with $g_0$, i.e.
$g_0(X,Y)=g_0(PX,PY)$ (for example, this is true for any canonical
almost product structure $P$ \cite{B2}). Clearly, $\mgot^-$ and
$\mgot ^+$ are orthogonal with respect to $g_0$, since for any $X
\in \mgot ^+, \; Y \in \mgot ^-$
$$g_0(X,Y)=g_0(P(X),P(Y))=g_0(X,-Y)=-g_0(X,Y).$$

Let us consider the action of the canonical almost product
structures on the $6$-symmetric space (\ref{space}). Here we use
notations of Corollary 2.

For $P_2(\theta)=\frac{1}{3}\theta+\theta ^2 +\frac{1}{3} \theta^3
+\theta ^4 +\frac{1}{3} \theta^5$ $\mgot ^-=\mgot _1 \cup \mgot
_2,$ $\mgot ^+=\mgot _3,$ therefore $\mgot _3 \bot \mgot _1$,
$\mgot _3 \bot \mgot _2$.

For $P_3(\theta)=\theta^3 $ $\mgot ^-=\mgot _1 \cup \mgot _3,$
$\mgot ^+=\mgot _2,$ thus $\mgot _2 \bot \mgot _1$. The statement
is proved. \hfill $\Box$

It can be deduced from Proposition 3 and  Proposition 5 that any
invariant Riemannian metric $g$ on (\ref{space}) is (up to a
positive coefficient) uniquely defined by the two positive numbers
$(s,t)$. It means that
\begin{equation}\label{metric}
g=g_0 |_{\mgot _1} +s g_0 |_{\mgot _2} + t g_0 |_{\mgot _3}.
\end{equation}

{\bf Definition 3.}  $(s,t)$ are called  the {\it characteristic
numbers} of the metric (\ref{metric}).

It should be pointed out that the canonical $f$-structures on the
homogeneous $\Phi$-space (\ref{space}) of the orders 4 and 6 are
metric $f$-structures with respect to all invariant Riemannian
metrics, which is proved by direct calculations.

Recall that in case of an arbitrary Riemannian metric $g$ the
Levi-Civita connection has its Nomizu function defined by the
formula (see \cite{KN})
\begin{equation}\label{alpha}
\alpha (X,Y)=\frac{1}{2}[X,Y] _{\mgot}+U(X,Y),
\end{equation}
where $X,Y \in \mgot$, the symmetric bilinear mapping $U$ is
determined by means of the formula
\begin{equation}\label{U}
2g(U(X,Y),Z)=g(X,[Z,Y]_{\mgot})+g([Z,X]_{\mgot},Y), \; X,Y,Z \in
\mgot .
\end{equation}

Suppose $g$ is an invariant Riemannian metric on the homogeneous
$\Phi$-spa\-ce (\ref{space}) with the characteristic numbers
$(s,t)$ $(s,t>0).$ The following statement is true.

{\bf Proposition\;6.} \begin{multline}\label{U1}
U(X,Y)=\frac{t-s}{2}([X_{\mgot _2},Y_{\mgot _3}]+[Y_{\mgot
_2},X_{\mgot _3}])+ \\ +\frac{t-1}{2s}([X_{\mgot _1},Y_{\mgot
_3}]+[Y_{\mgot _1},X_{\mgot _3}])+\frac{s-1}{2t}([X_{\mgot
_1},Y_{\mgot _2}]+[Y_{\mgot _1},X_{\mgot _2}]).
\end{multline}

{\bf Outline of the proof.} First we apply (\ref{metric}) and the
definition of $g_0$ to (\ref{U}). We take four matrices
$X=\{x_{i\, j}\}$, $Y=\{y_{i\, j}\}$, $Z=\{z_{i\, j}\}$ and
$U=\{u_{i\, j}\}$ and calculate the right-hand and left-hand side
of the equality obtained. After that we can represent it in the
form
\begin{equation}\label{zero}
c_{1\,2}z_{1\,2}+c_{1\,3} z_{1\,3}+\sum _{i=1} ^n c_{1\,i}
z_{1\,i}+\sum _{i=1} ^n c_{2\,i} z_{2\,i}+\sum _{i=1} ^n c_{3\,i}
z_{3\,i}=0,
\end{equation}
where $c_{1\,2},\;c_{1\,3},\;c_{1\,i},\;c_{2\,i},\;c_{3\,i}$
$(i=1,\dotsc,n)$ depend on elements of the matrices $X,Y$ and $U$.
As (\ref{zero}) holds for any $Z \in \mgot$, it follows in the
standard way that
\begin{equation}\label{c}
c_{1\,2}=c_{1\,3}=c_{1\,i}=c_{2\,i}=c_{3\,i}=0,\qquad
(i=1,\dotsc,n).
\end{equation}
Using (\ref{c}), we calculate $u_{i\,j}=u_{i\,j}(X,Y)$. To
conclude the proof, it remains to transform the formula for
$U(X,Y)$ into (\ref{U1}), which is quite simple.\hfill $\Box$

 In the notations of Section 2 we have the following statement.

{\bf Theorem\;2.} {\it Consider $SO(n)/SO(2) \times SO(n-3)$ as a
$4$-symmetric $\Phi$-space. Then the only canonical $f$-structure
$f_0$ on this space is}
\begin{itemize}
\item [1)] {\it a Killing $f$-structure iff
the characteristic numbers of a Riemannian metric are
$(1,\frac{4}{3})$;}
\item [2)] {\it a nearly K\"ahler $f$-structure iff
the characteristic numbers of a Riemannian metric are $(1,t), \;
t>0$;}
\item [3)] {\it a $G_1 f$-structure with respect to any invariant
Riemannian metric.}
\end{itemize}

{\bf Proof.} Application of (\ref{alpha}) to the definitions of
the classes ${\bf Kill \, f}$, ${\bf NKf}$ and ${\bf G_1 f}$
yields that
\begin{itemize}
\item [1)] $f \in {\bf Kill \, f}$ iff
$\frac{1}{2}[X,fX]_{\mgot}+U(X,fX)-f(U(X,X))=0$; \item [2)] $f \in
{\bf NKf}$ iff
$\frac{1}{2}[fX,f^2X]_{\mgot}+U(fX,f^2X)-f(U(fX,fX))=0$; \item
[3)] $f \in {\bf G_1 f}$ iff $f(2U(fX,f^2 X)-f(U(fX,fX))+f(U(f^2
X,f^2 X)))=0$.
\end{itemize}

The proof is straightforward. For example, it is known that $f_0$
is a nearly K\"ahler $f$-structure in the naturally reductive
case, which means that $\frac{1}{2}[f_0X,f_0^2 X]_{\mgot}=0$ for
any $X \in \mgot$ (see \cite{B6}). Making use of Proposition 4 and
Proposition 6, we obtain $U(f_0 X,f_0 X)\in \Ker f_0$ for any $X
\in \mgot$, $U(f_0 X,f_0^2 X)=0$ for any $X \in \mgot$ iff $s=1$.
Thus we have 2). Other statements are proved in the same manner.
\hfill $\Box$

 The similar technique is used to prove

{\bf Theorem 3.} {\it Consider  $SO(n)/SO(2)\times SO(n-3)$ as a
$6$-symmetric space. Let $(s,t)$ be the characteristic numbers of
an invariant Riemannian metric. Then}

\begin{itemize}
\item [1)]
 {\it $f_1$ is a Killing $f$-structure iff $s=1, \;
t=\frac{4}{3}$;\\
 $f_2$, $f_3$, $f_4$ do not belong to ${\bf Kill \, f}$ for any $s$
and $t$.}
\item [2)]
{\it $f_1$ is an $NKf$-structure iff s=1;\\
 $f_2$ and $f_3$ are $NKf$-structures for any $s$ and $t$;\\
 $f_4$ is not an $NKf$-structure for any $s$ and $t$.}
\item [3)] {\it $f_1$, $f_2$, $f_3$, $f_4$ are $G_1 f$-structures for any $s$ and
$t$.}
\end{itemize}

\end{document}